\documentclass[12pt,a4paper]{article}
\usepackage{latexsym,amsfonts,amsmath,graphics}
\usepackage{epsfig}

\newtheorem{theorem}{Theorem}

\newtheorem{definition}{Definition}

\newcommand{\bsa}{\boldsymbol{a}}

\newcommand{\bfn}{\mathbf{n}}

\newcommand{\bsx}{\boldsymbol{x}}
\newcommand{\bsh}{\boldsymbol{h}}

\newcommand{\bsz}{\boldsymbol{z}}

\newcommand{\bsy}{\boldsymbol{y}}

\newcommand{\cP}{{\cal P}}

\newcommand{\icomp}{\mathtt{i}}

\newcommand{\rd}{\,\mathrm{d}}
\newcommand{\NN}{\mathbb{N}}
\newcommand{\ZZ}{\mathbb{Z}}

\newcommand{\RR}{\mathbb{R}}

\newcommand{\FF}{\mathbb{F}}

\renewcommand{\pmod}[1]{\,(\bmod\,#1)}

\def\RR{\mathbb{R}}
\def\ZZ{\mathbb{Z}}
\def\NN{\mathbb{N}}

\def \F {\mathbb{F}}

\allowdisplaybreaks

\begin{document}

\title{Some highlights of Harald Niederreiter's work}
\author{Gerhard Larcher\thanks{G. Larcher is supported by the Austrian Science Fund (FWF): Project F5507-N26, which is a part of the Special Research Program "Quasi-Monte Carlo Methods: Theory and Applications".}, Friedrich Pillichshammer\thanks{F. Pillichshammer is supported by the Austrian Science Fund (FWF): Project F5509-N26, which is a part of the Special Research Program "Quasi-Monte Carlo Methods: Theory and Applications".}, Arne Winterhof\thanks{A. Winterhof is supported by the Austrian Science Fund (FWF): Project F5511-N26, which is a part of the Special Research Program "Quasi-Monte Carlo Methods: Theory and Applications".}\\ 
and Chaoping Xing\thanks{C. Xing is supported by Singapore Ministry of Education Tier 1 grant 2013-T1-002-033}}
\date{}
\maketitle

\begin{center}
Dedicated to our teacher, colleague and 
friend, Harald Niederreiter, on the occasion of his 70th birthday
\end{center}

\begin{abstract}
In this paper we give a short biography of Harald Niederreiter and we spotlight some cornerstones from his wide-ranging work. We focus on his results
on uniform distribution, algebraic curves, polynomials and quasi-Monte Carlo methods. In the flavor of Harald's work we also mention some applications
including numerical integration, coding theory and cryptography.
\end{abstract}

\section{A short biography }

Harald Niederreiter was born in Vienna in 1944 on June 7 and spent his childhood in Salzburg. 
In 1963 he returned to Vienna to study at the Department of Mathematics of the University of Vienna, 
where he finished his PhD thesis entitled ``Discrepancy in Compact Abelian Groups'' {\it sub auspiciis praesidentis rei publicae}
\footnote{The term ``Promotion sub auspiciis praesidentis rei publicae'' is the highest possible honor for course achievement at school and university in Austria.} 
under the supervision of Edmund Hlawka in 1969. From 1969 to 1978 he worked as scientist and professor in the USA at four different institutes:
Southern Illinois University, University of Illinois at Urbana-Champaign,
Institute for Advanced Study Princeton, and University of California at Los Angeles.
From 1978 to 1981 he was Chair of Pure Mathematics at the University of the West Indies in Kingston (Jamaica).
He returned to Austria and served as director of two institutes of the Austrian Academy of Sciences in Vienna,
until 1999 of the Institute for Information Processing and then of the Institute of Discrete Mathematics.
From 2001 to 2009 he was professor at the National University of Singapore.
Since 2009 he is located at the Johann Radon Institute for Computational and Applied Mathematics in Linz.
From 2010 to 2011 he was professor at the King Fahd University of Petroleum and Minerals in Dhahran (Saudi-Arabia).

Harald Niederreiter's research areas include
numerical analysis, pseudorandom number generation,
quasi-Monte Carlo methods,
cryptology, finite fields, applied algebra, algorithms, number theory, and
coding theory. He has published more than 350 research papers and several books including
\begin{itemize}
 \item  (with L. Kuipers) Uniform distribution of sequences. Pure and Applied Mathematics.
 Wiley-Interscience, New York-London-Sydney, 1974; reprint, Dover Publications, 2006.
 \item  (with R. Lidl) Finite fields. Encyclopaedia of Mathematics and its Applications, 20. Addison-Wesley Publishing Company, Advanced Book Program, Reading, MA, 1983; 2nd ed., 
 Cambridge University Press, 1997.
 \item (with R. Lidl) Introduction to Finite Fields and Their Applications, Cambridge University Press, 1986; revised ed., Cambridge University Press, 1994.
 \item  Random number generation and quasi-Monte Carlo methods. CBMS-NSF Regional Conference Series in Applied Mathematics, 63.
 Society for Industrial and Applied Mathematics (SIAM), 1992.
 \item  (with C.P. Xing) Rational points on curves over finite fields: theory and applications.
 London Mathematical Society Lecture Note Series, 285. Cambridge University Press, 2001.
 \item  (with C.P. Xing) Algebraic geometry in coding theory and cryptography.
 Princeton University Press, 2009.
 \end{itemize}
Furthermore he is editor or co-editor of 
\begin{itemize}
\item (with P.J.-S. Shiue) Monte Carlo and Quasi-Monte Carlo Methods in Scientific Computing, Springer-Verlag, 1995.
\item (with S.D. Cohen) Finite Fields and Applications, London Mathematical Society Lecture Note Series, 233, Cambridge University Press, 1996.
\item (with P. Hellekalek, G. Larcher und P. Zinterhof) Monte Carlo and Quasi-Monte Carlo Methods 1996, Springer-Verlag, 1998.
\item (with C. Ding und T. Helleseth) Sequences and Their Applications, Springer-Verlag, 1999.
\item (with J. Spanier) Monte Carlo and Quasi-Monte Carlo Methods 1998, Springer-Verlag, 2000.
\item (with D. Jungnickel) Finite Fields and Applications, Springer-Verlag, 2001.
\item (with K.-T. Fang and F.J. Hickernell) Monte Carlo and Quasi-Monte Carlo Methods 2000, Springer-Verlag, 2002.
\item Coding Theory and Cryptology, World Scientific Publishing, 2002.
\item Monte Carlo and Quasi-Monte Carlo Methods 2002, Springer-Verlag, 2004.
\item (with K. Feng und C.P. Xing) Coding, Cryptography and Combinatorics, Birkh\"auser-Verlag, 2004.
\item (with D. Talay) Monte Carlo and Quasi-Monte Carlo Methods 2004, Springer-Verlag, 2006.
\item (with A. Keller and S. Heinrich) Monte Carlo and Quasi-Monte Carlo Methods 2006, Springer-Verlag, 2008.
\item (with Y. Li, S. Ling, H. Wang, C.P. Xing and S. Zhang) Coding and Cryptology, World Scientific Publishing, 2008.
\item (with A. Ostafe, D. Panario and A. Winterhof) Algebraic Curves and Finite Fields: Cryptography and other applications, de Gruyter, 2014.
\item (with P. Kritzer, F. Pillichshammer and A. Winterhof) Uniform Distribution and Quasi-Monte Carlo Methods: Discrepancy, Integration and Applications, de Gruyter, 2014.
\end{itemize}

 Some important methods are named after him, such as the Niederreiter public-key cryptosystem,
 the Niederreiter factoring algorithm for polynomials over finite fields, and the Niederreiter and Niederreiter-Xing low-discrepancy sequences.

It follows an excerpt of his honors and awards:
\begin{itemize}
\item  full member of the Austrian Academy of Sciences
\item full member and former member of the presidium of the German Academy of Natural Sciences Leopoldina
\item Cardinal Innitzer Prize for Natural Sciences in Austria
\item invited speaker at ICM 1998 (Berlin) and ICIAM 2003 (Sydney)
\item Singapore National Science Award 2003
\item honorary member of the Austrian Mathematical Society 2012
\item Fellow of the American Mathematical Society 2013.
\end{itemize}

Niederreiter is also the initiator and, from 1994 to 2006, the co-chair of the first seven biennial {\it Monte Carlo and quasi-Monte Carlo meetings} which took place in:
\begin{itemize}
\item Las Vegas, USA (1994)
\item Salzburg, Austria (1996)
\item Claremont, USA (1998)
\item Hong Kong (2000)
\item Singapore (2002)
\item Juan-Les-Pins, France (2004)
\item Ulm, Germany (2006)
\item Montreal, Canada (2008)
\item Warsaw, Poland (2010)
\item Sydney, Australia (2012)
\item Leuven, Belgium (2014)
\end{itemize}
In 2006 Harald Niederreiter announced his wish to step down from the organizational role, and a Steering Committee was formed to ensure and oversee the continuation of the conference series.

\section{Uniform distribution theory and number theory}\label{secUDT}

When we scroll over the more than 350 scientific articles of Niederreiter which appeared in such renowned journals like "Mathematika", "Duke Mathematical Journal", "Bulletin of the AMS", or "Compositio Mathematica", we find that most of these papers have connections to topics from number theory or use techniques from number theory, and many of these articles deal with problems and solve open questions, or initiate a new field of research in the theory of uniform distribution of sequences. So also the later sections on Harald's work on coding theory, algebraic curves and function fields, on 
pseudorandom numbers, on finite fields, 
and on quasi-Monte Carlo methods in this overview in a certain sense will deal with number-theoretical aspects.

Let us give just one example: The analysis and the precise estimation of exponential sums $\sum\limits_{k=0}^{N-1}  {\rm e}^{2 \pi \icomp f(k)}$, or, 
in particular, of character sums plays an essential role in many different branches of mathematics and especially in number theory. So particularly it plays a basic role in 
many questions concerning uniform distribution of sequences, discrepancy theory, in quasi-Monte Carlo methods, in pseudorandom number analysis, in the theory of finite fields, 
and many more. In a variety of papers on exponential sums and their applications, Niederreiter has proven to be a leading expert in the analysis of exponential sums and has essentially developed various most important techniques.

In this section we just want to pick out some of the most impressive pieces of work of Niederreiter on some topics in number theory and in uniform distribution theory that will not be 
explicitly described in the subsequent sections.

In the first years after finishing his PhD thesis on "Discrepancy in Compact Abelian Groups" under the supervision of Edmund Hlawka, Niederreiter was concerned with basic questions from the theory of uniform distribution, from discrepancy theory and from metrical uniform distribution theory. We want to highlight three of the papers of this first phase.

In the paper ``An application of the Hilbert-Montgomery-Vaughan inequality to the metric theory of uniform distribution mod 1'' \cite{HilMont} which appeared in 1976 in the Journal of the London Mathematical Society, Niederreiter uses tools from the theory of bounded quadratic and bilinear forms, 
especially an inequality of Montgomery and Vaughan based on large sieve methods, to establish an analogue of Koksma's metric theorem for uniform distribution modulo one with respect to a general class of summation methods.

One of the most powerful tools for estimating discrepancy of sequences is the inequality of Koksma-Erd\H{o}s-Tur\'{a}n which bounds the discrepancy of a sequence by a weighted sum of the values of its Weyl sums. In the joint paper ``Berry-Esseen bounds and a theorem of Erd\H{o}s and Tur\'{a}n on uniform distribution mod 1'' \cite{NiedPhil} with Walter Philipp - this paper appeared in the Duke Mathematical Journal in 1973 - a much more general result about distances of functions is shown that contains the one-dimensional Koksma-Erd\H{o}s-Tur\'{a}n inequality as a special case. The given theorem is an analogue of the standard Berry-Esseen lemma for $\RR^s$.

One of the highlights in this period and of the work of Niederreiter in metric diophantine approximation theory certainly was the solution of a conjecture of Donald Knuth, together with Robert F. Tichy, in the paper ``Solution of a problem of Knuth on complete uniform distribution of sequences'' \cite{Tichy1} which appeared in Mathematika in 1985. It is shown there that for any sequence $(a_n)_{n \ge 1}$ of 
distinct positive integers, the sequence $(x^{a_n})_{n\ge 1}$ is completely uniformly distributed modulo one for almost all real numbers $x$ with $| x | > 1$. In the paper ``Metric theorems on uniform distribution and approximation theory'' \cite{Tichy2}, again in cooperation with Tichy, 
this result was even generalized to the following form: The sequence $(c x^{b_n})_{n\ge 1}$ is completely uniformly distributed modulo one for all $c \neq 0$ for almost all real numbers $x$ with $| x | > 1$ whenever $(b_n)_{n \ge 1}$ is any sequence of reals with $\inf b_n > -\infty $ and $\inf_{m \neq n} |b_n - b_m| > 0.$

In the analysis of the distribution properties of sequences and of point sets, especially of Kronecker sequences $$((\{n \alpha_1 \}, \ldots, \{n \alpha_s \}))_{n \ge 0}$$ and of lattice point sets $$\left(\left(\left\{n \frac{a_1}{N} \right\}, \ldots, \left\{n \frac{a_s}{N} \right\}\right)\right)_{n = 0, \ldots , N-1}$$ in the $s-$dimensional unit-cube, one often is led to questions from the theory of diophantine approximations, of the geometry of numbers or to questions concerning continued fraction expansions. A famous still open problem in the theory of continued fractions is the following conjecture of Zaremba:

{\it There is a constant $c$ such that for every integer $N \ge 2$ there exists an integer $a$ with $1 \le a \le N$ and with $\gcd(a, N) = 1$ such that all continued fraction coefficients of $\frac{a}{N}$ are bounded by $c$. Indeed it is conjectured that $c = 5$ satisfies this property.}

Niederreiter in ``Dyadic fractions with small partial quotients'' \cite{Mh} proved that this result is true even with $c=3$ if $N$ is a power of 2. He proves the conjecture of Zaremba also for $N$ equal to powers of 3 and equal to powers of 5. Only quite recently it was shown by Bourgain and Kontorovich that Zaremba's conjecture holds for almost all choices of $N$.

From Niederreiter's result now for example it can be deduced that for all $N = 2^m$ there exists an integer $a$ such that the lattice point set $$\left(\left(\left\{n \frac{1}{2^m} \right\},  \left\{n \frac{a}{2^m} \right\}\right)\right)_{n = 0, \ldots , 2^m-1}$$ has discrepancy $D_N \le c' \frac{\log N}{N}$, i.e., has best possible order of discrepancy.

The investigation of certain types of digital $(t,m,s)$-nets and of digital $({\bf T},s)$-sequences (see also Section~\ref{secQMC} of this article) in analogy leads to questions concerning non-archimedean diophantine approximation and to questions concerning continued fraction expansions of formal Laurent series. Such questions were analyzed, for example, in the papers \cite{ActaArith1993,LarNiedTransact1995,MoMa1987,SciMathHungar1995}.

In an impressive series of papers together with Igor Shparlinski, powerful methods for the estimation of exponential sums with nonlinear recurring sequences were developed by Niederreiter,
see also Section~\ref{poly} below. In the paper ``On the distribution of power residues and primitive elements in some nonlinear recurring sequences'' \cite{BullLonSoc2003} which appeared in the Bulletin of the London Mathematical Society in 2003 it is shown that these methods also can be applied to estimating sums of multiplicative characters as well.
As a consequence, results are obtained in this paper about the distribution of power residues and of primitive elements in such sequences.

So consider a sequence of elements $u_0, u_1, \ldots, u_{N-1}$ of the finite field $\FF_q$ obtained by the recurrence relation $$u_{n+1}= a u_n^{-1} +b,$$
where we set $u_{n+1} = b$ if $u_n = 0$. For a divisor $s$ of $q-1$ let $R_s(N)$ be the number of $s$-power residues (i.e., the number of $w \in \FF_q$ such that there are $z \in \FF_q$ with $z^s = w$) among $u_0, u_1, \ldots, u_{N-1}$. 
Then
$$\left| R_s(N) - \frac{N}{s}\right| < (2.2)  N^{\frac{1}{2}} q^{\frac{1}{4}}$$
for $1 \le N \le t$, where $t$ is the least period of the recurring sequence. The case of general nonlinear recurrence sequences was studied later in \cite{niwiaa}.

Also in the present Harald Niederreiter is still an utmost creative and productive researcher in the field of number theory and uniform distribution of sequences. We want to confirm this fact by giving two last quite recent examples of his impressive work in these fields:

In the joint paper ``On the Gowers norm of pseudorandom binary sequences'' \cite{BullAusMathSoc2009} with Jo\"el Rivat, the modern concepts of Christian Mauduit and Andr\'{a}s S\'{a}rk\"ozy 
concerning new measures for pseudorandomness and of William T. Gowers in combinatorial and additive number theory were brought together and the Gowers norm for periodic binary sequences is studied. 
A certain relation between the Gowers norm of a binary function $f$ defined on the integers modulo $N$ and a certain correlation measure for the sequence $(f(n))_{n \ge 1}$
introduced in \cite{masa} is shown.

A quite new and challenging trend in the theory of uniform distribution of sequences is the investigation of the distribution of hybrid sequences. A hybrid sequence is defined as follows:
take an $s$-dimensional sequence $(\bsx_n)_{n \ge 0}$ of a certain type and a $t$-dimensional sequence $(\bsy_n)_{n \ge 0}$ of another type and combine them to an $(s+t)$-dimensional {\em hybrid sequence}, i.e., with some abuse of notation, $$(\bsz_n)_{n \ge 0}  : = ((\bsx_n,\bsy_n))_{n \ge 0}.$$
Well known examples of such sequences are Halton-Kronecker sequences (generated by combining Halton sequences with Kronecker sequences) or Halton-Niederreiter sequences 
(a combination of digital $(t,s)$-sequences or of digital $({\bf T},s)$-sequences in different bases). The investigation of these sequences again leads to challenging problems in number theory. For example with his papers \cite{AA2009,UDT2010,MM2010,Debrecen2011,AA2012} Niederreiter
essentially influences the direction of research in this topic.

\section{Algebraic curves, function fields and applications}\label{sec:FF}
The study of algebraic curves over finite fields can be traced back to Carl Friedrich Gau{\ss} who studied equations over finite fields. However, the real beginning of this topic was the proof of the Riemann hypothesis for algebraic curves over finite fields by 
Andr\'{e} Weil in the 1940s. This topic has attracted great attention of the researchers again since the 1980s due to the discovery of algebraic geometry codes  by Valerii D. Goppa. This application of algebraic curves over finite fields, and especially of those with many rational points, created a much stronger interest in the area and attracted new groups of researchers such as coding theorists and algorithmically inclined mathematicians. Nowadays, algebraic curves over finite fields is a flourishing subject  which produces exciting research and is immensely relevant for applications.

Harald Niederreiter started this topic from applications first. In the late 1980s, he found an elegant construction of $(t,m,s)$-nets and  $(t,s)$-sequences 
(see Section \ref{secQMC}). Then he realized that the construction can be generalized to global function fields \cite{NX96,NX96a}. From this point, Harald Niederreiter has extensively investigated algebraic curves over finite fields with many rational points and their applications.
Algebraic curves over finite fields can be described in an equivalent algebraic language, i.e, global function fields over finite fields. For many of the applications, people are interested in algebraic curves over finite fields with many rational points, or equivalently, global function fields over finite fields with many rational places. Since the global function field language was usually used by Harald Niederreiter, we adopt this language from now onwards  in this section.

Let $\F_q$ denote the finite field of $q$ elements. An extension $F$ of $\F_q$ is called an algebraic function field of one variable over $\F_q$ if there exists an element $x$ of $F$ that is transcendental over $\F_q$ such that $F$ is a finite extension over the rational function field $\F_q(x)$. We usually denote by $F/\F_q$ a global function field with the full constant field $\F_q$, i.e., all  elements in $F\setminus\F_q$ are transcendental over $\F_q$.  A place $P$ of $F$ is called {\it rational} if its residue field $F_P$ is isomorphic to the ground field $\F_q$. For many applications in coding theory, cryptography and low-discrepancy sequences, people are interested in those function fields with many rational places. On the other hand, the number of rational places of a function field over $\F_q$ is constrained by an important invariant of $F$, called genus. If we use $g(F)$ and $N(F)$ to denote genus and number of rational places of $F/\F_q$, the well-known Hasse-Weil bound says that
\begin{equation}\label{eqn:HW}
|N(F)-q-1|\le 2g(F)\sqrt{q}.\end{equation}
The above bound implies that the number of rational places cannot be too big if we fix the genus of a function field. Now the game becomes to find the maximal number of rational places that a global function field over $\F_q$ of genus $g$ could have. We usually denote by $N_q(g)$ this quantity, i.e,
$N_q(g)=\max\{N(F):\; F/\F_q \ \mbox{has genus} \ g\}$. Apparently, it follows from the Hasse-Weil bound that
\begin{equation}\label{eqn:HW1}
|N_q(g)-q-1|\le 2g\sqrt{q}\end{equation}
for any prime power $q$ and nonnegative integer $g$. For given $q$ and $g$, determining the exact value of $N_q(g)$  is a major problem in the study of global function fields. 
In general it is very difficult to determine the exact value of $N_q(g)$. Instead, it is sufficient to find reasonable lower bounds for most applications. Finding lower bounds on $N_q(g)\ge N$ is through either explicit construction or showing existence of global function fields of genus $g$ with at least $N$ rational places. Investigation of this problem involves several subjects such as algebraic number theory and algebraic geometry and even coding theory. The method that Harald Niederreiter employed is class field theory in algebraic number theory. He found many record function fields through class field theory, i.e., global function fields with best-known number of rational places. Some of these record function fields are listed below (see \cite{NX96a,N-X2,N-X3,N-X4,N-X5,N-X6,N-X7,NX98,X-N}).
\begin{center}
\begin{tabular}{|c|c|c|c|c|c|c|c|c|c|c|c|}\hline
$(q,g)$ & $(2,23)$ & $(2,25)$& $(2,29)$ & $(2,31)$ & $(2,34)$ &$(2,36)$ & $(2,49)$ & $(3,6)$ & $(3,7)$ \\ \hline
$N_q(g)$ & $22$ & $24^*$ & $25$ &          $27$ &      $27$ & $30$ &      $36$ & $14^* $  &   $16^*$ \\ \hline
\end{tabular}
\end{center}
The entries with asterisk in the above table are the exact values of $N_q(g)$, while the entries without asterisk in the above table  are lower bounds on $N_q(g)$.

For a fixed prime power $q$, to measure how $N_q(g)$ behaves while $g$ tends to infinity, we define the following asymptotic quantity
 \begin{equation}\label{eqn:A(q)}
A(q):=\limsup_{g\rightarrow\infty}\frac{N_q(g)}g.\end{equation}
It is immediate from the Hasse-Weil bound that $A(q)\le 2\sqrt{q}$. Sergei G. Vl\u{a}du\c{t} and Vladimir G. Drinfeld refined this bound to $A(q)\le \sqrt{q}-1$. Yasutaka Ihara first showed that $A(q)\ge \sqrt{q}-1$ if $q$ is a square. Thus, the problem of determining $A(q)$ is completely solved for squares $q$. 
It still remains to determine $A(q)$ for nonsquare $q$. Like the case of $N_q(g)$, finding the exact value of $A(q)$ for nonsquare $q$ is very difficult in this topic. Although people have tried very hard, so far $A(q)$ has not been determined for any single nonsquare $q$. In particular, 
if $q$ is a prime, it is a great challenge to determine or find a reasonable lower bound on $A(q)$.

What Harald Niederreiter did on this problem was finding a new bound on $A(2)$ and an improvement on $A(q^m)$ for odd $m$. More precisely, he proved the following result in \cite{NX98d,NX99}.
\begin{theorem}
One has $A(2)\ge \frac{81}{317}=0.2555....$
\end{theorem}

\begin{theorem} One has the following bounds:
\begin{itemize}
\item[{\rm (i)}] 
 {\it If q is an odd prime power and} $m\ge 3$ {\it is an integer, then}
$$A(q^{m})\ge\frac{2q+2}{\lceil 2(2q+3)^{1/2}\rceil+1}.$$
\item[{\rm (iI)}] {\it If} $q\ge 8$ {\it is a power of} {\rm 2} {\it and} $m\ge 3$ {\it is an odd
integer, then}
$$A(q^{m})\ge\frac{q+1}{\lceil 2(2q+2)^{1/2}\rceil+2}.$$
\end{itemize}
\end{theorem}

Harald Niederreiter has been also working on applications of algebraic curves over finite fields. These applications include low-discrepancy sequences, coding theory and cryptography, etc. For the details on application of algebraic curves over finite fields to low-discrepancy sequences, we refer to Section \ref{secQMC}.

For applications of algebraic curves over finite fields to coding theory, Harald Niederreiter's contribution was the discovery of several new codes via the theory of algebraic curves over finite fields. Some of the new codes discovered by Harald Niederreiter are listed below (see \cite{DNX00}). In the following table, $[n,k,d]_q$ means a $q$-ary code of length $n$, dimension $k$ and minimum distance $d$. 
\begin{center}
\begin{tabular}{|c|c|c|c|c|c|}\hline
$[108,25,44]_4$ & $[108,26,43]_4$ & $[113,27,45]_4$ & $[130,29,53]_4$ & $[27,11,13]_8$ & $[30,7,19]_8$\\ \hline
 $[30,8,18]_8$ & $[30,9,17]_8$ & $[36,7,23]_8$   & $[36,8,22]_8$  & $[36,9,21]_8$  & $[36,10,20]_8$   \\ \hline
\end{tabular}
\end{center}

Harald Niederreiter has also done some significant work on asymptotic results of coding theory and cryptography as well via  algebraic curves over finite fields.

\section{Polynomials over finite fields and applications}\label{poly}

Now we describe some of Harald Niederreiter's results on polynomials over finite fields and applications.
 We start with complete mappings and check digit systems.
 
  Let $\F_q$ be the finite field of $q>2$ elements and 
$f(X)\in\F_q[X]$ a permutation polynomial over $\F_q$. 
We call $f(X)$ a {\em complete mapping}
if 
$f(X)+X$ 
is also a permutation polynomial. Existence results on complete mappings and their application to check digit systems were discussed in
\cite{niro,shwi}. 

It is easy to see that $f(X)=aX$ is a complete mapping whenever $a\not\in \{-1,0\}$. 

Complete mappings are pertinent to the construction of orthogonal Latin squares, see \cite{ma}, which can be used to design some agricultural experiments. However, here we will
describe another application of complete mappings, namely, check digit systems. 

A {\em check digit system} (defined with one permutation polynomial over $\F_q$) 
consists of a permutation polynomial $f(X)\in \F_q[X]$ and a control symbol $c\in \F_q$ such that
each word $a_1,\ldots,a_{s-1}\in \F_q^{s-1}$ of length $s-1$ is extended by a check digit $a_s\in \F_q$ such 
that 
$$\sum_{i=0}^{s-1} f^{(i)}(a_{i+1})=c,$$
where $f^{(i)}$ is recursively defined by $f^{(0)}(X)=X$ and $f^{(i)}(X)=f(f^{(i-1)}(X))$ for $i=1,2,\ldots$

An example of a check digit system is 
the international standard book number (ISBN-10)
which consists of a string of $10$ digits $x_1-x_2x_3x_4x_5x_6-x_7x_8x_9-x_{10}$. 
 The first digit~$x_1$ characterizes the language group, $x_2x_3x_4$   is the number
of the publisher,  $x_5x_6x_7x_8x_9$ is the actual book number, and $x_{10}$ is a check digit. A correct ISBN satisfies 
$$x_1+2x_2+3x_3+4x_4+5x_5+6x_6+7x_7+8x_8+9x_9+10x_{10}=0\in \F_{11}.$$
With the variable transformation $a_i= x_{2^{i-1}\bmod 11}$ we get a check digit system defined with one permutation $f(X)=2X$. 
Note that $f(X)=2X$ and $-f(X)=9X$ are both complete mappings of $\F_{11}$.

For example the ISBN-10 of the monograph on finite fields by Lidl and Niederreiter \cite{lini} is
$0-521-39231-4$.

Since $f(X)$ is a permutation polynomial, such a system detects all single errors $a\mapsto b$.
Moreover it detects all  
\begin{itemize}
 \item neighbor transpositions $ab\mapsto ba$ if $-f(X)$ is a complete mapping;
 \item twin errors $aa\mapsto bb$ if $f(X)$ is a complete mapping.
\end{itemize}

Niederreiter and Karl H. Robinson \cite{niro} found several nontrivial classes of complete mappings and proved in particular a generalization of the following result:
\begin{theorem}
  Let $q$ be odd. Then $f_b(X)=X^{(q+1)/2}+bX$ is 
  a complete mapping of $\F_q$ if and only if $b^2-1$ and $b^2+2b$ are both squares of nonzero elements of $\F_q$.
  The number of $b$ such that $f_b(X)$ is a complete mapping is $\frac{q}{4}+O(q^{1/2})$.  
\end{theorem}
For a survey on (generalizations of) complete mappings and some applications we refer to \cite{wi}.\\

 Harald Niederreiter also invented a {\em deterministic algorithm} based on linear algebra for {\em factoring a univariate polynomial} $f(X)$ over $\F_q$ which is efficient for small characteristic,
see \cite{ni93} for the initial article and \cite{gapa} for a survey on factorization. The key step is to find a polynomial $h(X)$ which satisfies the differential equation
$$f^q(h/f)^{(q-1)}+h^q=0,$$
where $g^{(k)}$ denotes the $k$th Hasse-Teichm\"uller derivative.
Then $\gcd(f,h)$ is a nontrivial factor of $f$.\\

Harald Niederreiter did not only contribute to cryptography via the above mentioned public-key cryptosystem named after him, but also in many other ways. For example he proved several results on the 
interpolation of the discrete logarithm \cite{ni90,niwi} showing that there is no low degree polynomial $f(X)\in \F_q[X]$ which coincides with the discrete logarithm on many values, 
that is for prime $q$,
$f(g^x)=x$ for many $x$, where $g$ is a primitive element of $\F_q$. Hence, the discrete logarithm problem is not attackable via simple interpolation which is necessary for the security of
discrete logarithm based cryptosystems as the Diffie-Hellman key exchange. \\

Finally, he introduced and studied {\em nonlinear pseudorandom number generators}, i.e., sequences over $\F_q$ of the form
$$u_{n+1}=f(u_n),\quad n=0,1,\ldots$$
for some initial value $u_0\in \F_q$ and a polynomial $f(X)\in \F_q[X]$ of degree at least~$2$. These sequences are attractive alternatives to linear pseudorandom number generators 
which are not suitable for all 
applications. For example, linear generators are highly predictable and not suitable in cryptography.

As mentioned before in joint work with Igor Shparlinski \cite{nish99,nish00} Niederreiter found a way to prove nontrivial estimates on certain character sums which in the simplest case are of the form
$$\sum_{n=0}^{N-1} \chi(f(u_n))$$
where $\chi$ is any nontrivial additive character of $\F_q$. For such character sums the standard method for estimating incomplete character sums by reducing them to complete ones 
and then applying the Weil bound does not work. The method and result of \cite{nish99} was later slightly improved in \cite{niwi08}. In particular if $f(X)=aX^{q-2}+b$, 
i.e.\ $f(c)=ac^{-1}+b$ if $c\ne 0$,  this method 
yields strong bounds on the exponential sums and leads to very good discrepancy bounds for corresponding sequences in the unit interval. 
For a survey on nonlinear recurrence sequences see \cite{wiseta}.

\section{Quasi-Monte Carlo methods}\label{secQMC}

The quasi-Monte Carlo method has its roots in the theory of uniform distribution modulo 1 (see Section~\ref{secUDT}) and is nowadays a powerful tool in computational mathematics, in particular for the numerical integration of very high dimensional functions, with many applications to practical problems coming from biology,
computer graphics, mathematical finance, statistics, etc. Here the integral of a function $f:[0,1]^s \rightarrow \RR$ is approximated by a quasi-Monte Carlo (QMC) rule which computes the arithmetic mean of function values over a finite set of sample nodes, i.e., $$\int_{[0,1]^s} f(\bsx) \rd \bsx \approx \frac{1}{N}\sum_{n=0}^{N-1} f(\bsx_n)$$ with fixed $\bsx_0,\ldots,\bsx_{N-1} \in [0,1)^s$. QMC rules can be viewed as deterministic versions of Monte Carlo rules. The fundamental error estimate for QMC rules is the Koksma-Hlawka inequality which bounds the absolute integration error as $$\left|\int_{[0,1]^s} f(\bsx) \rd \bsx - \frac{1}{N}\sum_{n=0}^{N-1} f(\bsx_n)\right| \le V(f) D_N^{\ast}(\bsx_0,\ldots,\bsx_{N-1}),$$ where $V(f)$ is the variation of $f$ in the sense of Hardy and Krause and where $D_N^{\ast}$ is the star discrepancy of the underlying sample nodes, see \cite{kuinie}.

In the mid of the 1970s Harald Niederreiter started to investigate QMC methods. His first pioneering work was the paper ``Quasi-Monte Carlo methods and pseudo-random numbers'' \cite{nie78} published in the Bulletin of the American Mathematical Society in 1978. Today this paper can be seen as the first systematic survey about the theoretical foundations of QMC dealing with Koksma-Hlawka type inequalities and with constructions of point sets for QMC rules such as Halton's sequence, Sobol's construction of $P_\tau$ nets and $LP_\tau$ sequences, and good lattice points in the sense of Korobov and Hlawka.

The quintessence of the Koksma-Hlawka inequality is that good QMC rules should be based on sample nodes with low discrepancy, informally often called {\it low-discrepancy point sets}. Today there are two main streams of constructing low-discrepancy point sets. Both constructions are intimately connected with the name Niederreiter who contributed pioneering works to these topics. The first construction is the concept of lattice point sets and the second one is the concept of $(t,m,s)$-nets and $(t,s)$-sequences in a base $b$. \\

An {\it $N$-element lattice point set} (cf. Section~\ref{secUDT}) is based on an $s$-dimensional lattice point $\bsa=(a_1,\ldots,a_s)$. The $n$th element of such a lattice point set is then given as $$\bsx_n=\left\{\frac{n}{N} \bsa\right\}\ \ \ \mbox{ for }\ n=0,1,\ldots,N-1,$$ where the fractional part function $\{\cdot \}$ is applied component-wise. QMC rules which are based on good lattice point sets are called the {\it method of good lattice points} or {\it lattice rules} and belong nowadays to the most popular QMC rules in practical applications. Niederreiter analyzed distribution properties and showed the existence of good lattice point sets with low discrepancy. The full power of lattice rules however lies in the integration of smooth one-periodic functions. One reason for this is the following relation: for $\bsh \in \ZZ^s$ $$\frac{1}{N} \sum_{n=0}^{N-1} \exp\left(2 \pi \icomp \frac{n}{N} \bsa \cdot \bsh\right)=\left\{
\begin{array}{ll}
1 & \mbox{ if }\ \bsa \cdot \bsh \equiv 0 \pmod{N},\\
0 & \mbox{ if }\ \bsa \cdot \bsh \not\equiv 0 \pmod{N},
\end{array}
\right.$$
where $\cdot$ denotes the usual inner product. Niederreiter studied the worst-case error $P_\alpha$ for the integration of functions $f$ which can be represented by absolutely convergent Fourier series whose Fourier coefficients $\widehat{f}(\bsh)$ tend to zero as $\bsh$ moves away from the origin at a prescribed rate which is determined by the parameter $\alpha$. His most important contributions to the theory of good lattice point sets are summarized in Chapter~ 5 of his book ``Random number generation and quasi-Monte Carlo methods'' \cite{niesiam} which appeared in 1992. Niederreiter's most recent contributions to the theory of lattice point sets deal with the existence and construction of so-called extensible lattice point sets which have the property that the number of points in the node set may be increased while retaining the existing points (see \cite{hicknie,niepill}). \\

The theory of {\it $(t,m,s)$-nets and $(t,s)$-sequences} was initiated by Niederreiter in his seminal paper ``Point sets and sequences with small discrepancy'' \cite{nie87} published in the Monatshefte f\"{u}r Mathematik in 1987. The basic idea of these concepts is that if a point set has good equidistribution properties with respect to a reasonable (finite) set of test sets, then the point set already has low star discrepancy at all. The definition of a $(t,m,s)$-net in base $b$ can be stated as follows:

\begin{definition}[Niederreiter, 1987]
Let $s,b,m,t$ be integers satisfying $s \ge 1$, $b \ge 2$ and $0 \le t \le m$. A set $\cP$ consisting of $b^m$ elements in $[0,1)^s$ is said to be a $(t,m,s)$-net in base $b$ if every so-called elementary interval of the form $$\prod_{j=1}^s \left[\frac{a_j}{b^{d_j}}, \frac{a_j+1}{b^{d_j}}\right)$$ of volume $b^{t-m}$ with $d_j \in \NN_0$ and $a_j \in \{0,1,\ldots,b^{d_j}-1\}$ for $j=1,2,\ldots,s$, contains exactly $b^t$ elements of $\cP$.
\end{definition}

A $(t,s)$-sequence in base $b$ is an infinite version of $(t,m,s)$-nets.

\begin{definition}[Niederreiter, 1987]
Let $s,b,t$ be integers satisfying $s \ge 1$, $b \ge 2$ and $t \ge 0$. An infinite sequence $(\bsx_n)_{n \ge 0}$ of points in $[0,1)^s$ is said to be a $(t,s)$-sequence in base $b$ if, for all integers $k \ge 0$ and $m >t$, the point set consisting of the $\bsx_n$ with $k b^m \le n < (k+1)b^m$ is a $(t,m,s)$-net in base $b$.
\end{definition}

In his work \cite{nie87} Niederreiter presented a comprehensive theory of $(t,m,s)$-nets and $(t,s)$-sequences including discrepancy estimates, existence results and connections to other mathematical disciplines such as, e.g., Combinatorics. The fundamental discrepancy estimate for a $(t,m,s)$-net $\cP$ in base $b$ states that $$D_N^{\ast}(\cP) \le c_{s,b} b^t \frac{(\log N)^{s-1}}{N} + O_{s,b}\left(b^t \frac{(\log N)^{s-2}}{N}\right)$$ where $N=b^m$ and where $c_{s,b}>0$ is independent of $m$ and $t$. This estimate justifies the definition of $(t,m,s)$-nets since it means that for sufficiently small $t$ one can achieve a star discrepancy of order of magnitude $O((\log N)^{s-1}/N)$. Many people from discrepancy theory conjecture that this is the best convergence rate at all which can be achieved for the star discrepancy of $N$-element point sets in dimension $s$. For infinite $(t,s)$-sequences in base $b$ one can achieve a star discrepancy of order of magnitude $O((\log N)^s/N)$ which again is widely believed to be the best rate at all for the star discrepancy of infinite sequences in dimension $s$.

Most constructions of $(t,m,s)$-nets and $(t,s)$-sequences rely on the digital method which was introduced by Niederreiter, also in \cite{nie87}. In the case of $(t,m,s)$-nets this construction requires $m \times m$ matrices $C_1,C_2,\ldots,C_s$ over a commutative ring $R$ with identity and $|R|=b$ and, in a simplified form, a bijection $\psi$ from the set of $b$-adic digits $\mathcal{Z}_b=\{0,1,\ldots,b-1\}$ onto $R$. For $n=0,1,\ldots,b^m-1$, 
let $n=n_0 +n_1 b+\cdots+n_{m-1} b^{m-1}$ with all $n_r \in \mathcal{Z}_b$. Then, for $j=1,2,\ldots,s$, multiply the matrix $C_j$ with the vector $\bfn=(\psi(n_0),\psi(n_1),\ldots,\psi(n_{m-1}))^\top$ whose components belong to $R$, $$C_j \bfn= (y_{n,j,1},y_{n,j,2},\ldots,y_{n,j,m})^\top,\ \ \ \mbox{ with all }\ y_r \in R$$ and set $\bsx_n=(x_{n,1},x_{n,2},\ldots,x_{n,s})$, where $$\bsx_{n,j}=\frac{\psi^{-1}(y_{n,j,1})}{b}+\frac{\psi^{-1}(y_{n,j,2})}{b^2}+\cdots +\frac{\psi^{-1}(y_{n,j,m})}{b^m}.$$ The point set $\{\bsx_0,\bsx_1,\ldots,\bsx_{b^m-1}\}$ constructed this way is a $b^m$-element point set in $[0,1)^s$ and it is therefore a $(t,m,s)$-net in base $b$ for some $t \in \{0,1,\ldots,m\}$ which is called a digital $(t,m,s)$-net over $R$. In the case of $(t,s)$-sequences the only difference is that one uses $\infty \times \infty$ matrices.

The so-called quality parameter $t$ depends only on the chosen matrices $C_1,C_2,\ldots,C_s$. Of course $t$ should be as small as possible, in the optimal case $t=0$. If the base $b$ is a prime power, then one chooses for $R$ the finite field $\FF_b$ of order $b$ which makes life a bit easier and is therefore the most studied case. Then $t$ is determined by some linear independence property of the row-vectors of the generating matrices $C_1,C_2,\ldots,C_s$ which provides the link of digital nets and sequences to the theory of finite fields and linear algebra over finite fields.

Niederreiter developed several constructions of generating matrices which lead to good, often even optimal small $t$-values. One important construction results in the now so-called {\it Niederreiter sequences} and is based on polynomial arithmetic over finite fields and the formal Laurent series expansion of certain rational functions over $\FF_b$ whose Laurent coefficients are used to fill the generating matrices. If $s \le b$ this leads to an explicit construction of $(0,s)$-sequences in base $b$ which in turn implies, for $s \le b+1$, an explicit construction of a $(0,m,s)$-net in base $b$ for every $m \ge 2$. It is known that the conditions $s \le b$ for sequences and $s \le b+1$ for nets, respectively, are even necessary to achieve a quality parameter equal to zero. Niederreiter sequences and slight generalizations thereof recover and unify the existing constructions due to Il'ya M. Sobol' and Henri Faure.

An important subclass of $(t,m,s)$-nets which was introduced by Niederreiter in the paper ``Low-discrepancy point sets obtained by digital constructions over finite fields'' \cite{nie92} is provided by the concept of what we call today {\it polynomial lattice point sets}. This name has its origin in a close relation to ordinary lattice point sets. In fact, the research on polynomial lattice point sets and on ordinary lattice point sets often follows two parallel tracks and bears a lot of similarities (but there are also differences).

Niederreiter's early work on $(t,m,s)$-nets and $(t,s)$-sequences is well summarized in Chapter~4 of his already mentioned book ``Random number generation and quasi-Monte Carlo methods'' \cite{niesiam} which appeared in 1992. Since its appearance, this book is {\it the} reference book for $(t,m,s)$-nets and $(t,s)$-sequences especially and for QMC and random number generation in general.

A disadvantage of Niederreiter sequences in dimension $s$ is that they can only achieve a $t$-value of order $O(s \log s)$ as $s$ tends to infinity. This disadvantage was overcome in the next cornerstone of Niederreiter's work in QMC, the constructions of Niederreiter-Xing sequences. In a series of papers \cite{NX95,NX96,NX96a,NX98,NX02} starting in 1995 and based on methods from Algebraic Geometry, Niederreiter developed in collaboration with Chaoping Xing constructions of generating matrices which achieve the currently best known quality parameters of order $O(s)$ for growing dimensions $s$. This order is known to be best possible. An introduction into this subject and an overview can be found in the book ``Rational points on curves over finite fields'' \cite{NX} published by Niederreiter and Xing in 2001.

In 2001 Niederreiter developed together with Gottlieb Pirsic \cite{NiePir} a {\it duality theory} for digital nets. The basic idea is that the construction of digital $(t,m,s)$-nets over $\FF_b$ can be reduced to the construction of certain $\FF_b$-linear subspaces of $\FF_b^{sm}$. Using the standard inner product in $\FF_b^{sm}$ one can define and study the dual linear subspace. If one defines a special weight on $\FF_b^{sm}$, the so-called {\it Niederreiter-Rosenbloom-Tsfasman weight}, then the $t$-parameter of a digital net is closely related to the weight of the corresponding dual linear subspace. This point of view gives new possibilities for the construction of digital nets, as, for example, cyclic nets or hyperplane nets, and it provides a connection to the theory of linear codes. Later, in 2009, Niederreiter and Josef Dick \cite{dienie} extended the duality theory for digital nets also to digital sequences which became a convenient framework for the description of many constructions such as, for example, the ones of Niederreiter and Xing or of Niederreiter and 
Ferruh \"Ozbudak \cite{NieOe02,NieOe04} (see also \cite{DP}).

Digital nets also have a close connection to other discrete objects such as orthogonal Latin squares or ordered orthogonal arrays. Also these relations were subject of Niederreiter's research.\\

Harald Niederreiter's contributions to the theory of QMC are groundbreaking. He opened new doors and developed comprehensive theories of lattice rules and of $(t,m,s)$-nets and $(t,s)$-sequences with many new ideas and facets. Today Niederreiter's work forms one of the essential pillars of QMC integration.

\end{document}